\documentclass[11pt, a4paper]{article}
\usepackage{amsmath}
\usepackage{amssymb,amsthm}
\usepackage[english]{babel}
\usepackage[utf8]{inputenc}
\usepackage[T1]{fontenc}
\usepackage{lmodern}
\usepackage{color}
\usepackage[total={15cm,25.2cm}, top=3cm, left=3cm, includefoot]{geometry}
\usepackage[ruled,linesnumbered]{algorithm2e}
\usepackage{hyperref}

\newcommand*\mystrut[1]{\vrule width0pt height0pt depth#1\relax}

\title{S-SPADE Done Right:\\ Detailed Study of the Sparse Audio Declipper Algorithms
	}
\author{Pavel Záviška, Ondřej Mokrý, Pavel Rajmic%
	\thanks{All the authors are with  Brno University of Technology, Czech Republic.
					Contact: {rajmic@vutbr.cz}.
					Authors are thankful to Jakub Kůdela for discussing the Augmented Lagrangian for complex variables.
}%
}
\date{\today}

\theoremstyle{plain}
\def\RR{\mathbb{R}}

\def\CC{\mathbb{C}}

\theoremstyle{plain}
\newtheorem{Rem}{Remark}

\def\x{\vect{x}}
\def\u{\vect{u}}

\def\y{\vect{y}}
\def\z{\vect{z}}
\def\w{\vect{w}}

\def\c{\vect{c}}

\def\rr{\vect{r}}
\def\s{\vect{s}}
\def\Mr{M_\text{R}}
\def\Mh{M_\text{H}}
\def\Ml{M_\text{L}}

\def\RA{\mathcal{R}(A)}

\newcommand{\Id}{\mathit{Id}}

\newcommand{\norm}[1]{\left\|#1\right\|}

\newcommand{\adjoint}[1]{#1^*}

\newcommand{\vect}[1]{\mathbf{#1}} %vector
\newcommand{\argmin}{\mathop{\operatorname{argmin}}}

\newcommand{\ig}{\iota_{\Gamma(\y)}(\x)}
\newcommand{\il}{\iota_{\ell_0 \leq k}(\z)}
\newcommand{\modulo}{\mathop{\operatorname{mod}}}

\begin{document}
\maketitle

\noindent
\emph{Abstract:}
This technical report shows and discusses in detail how Sparse Audio Declipper (SPADE) algorithms are derived from the signal model using the ADMM approach.
The analysis version (A-SPADE) of Kiti\'c et.\,al.\ (LVA/ICA\,2015) is derived and justified.
The synthesis version (S-SPADE) of the same research team is shown to solve a~different optimization task than intended.
This issue is corrected in this report, leading to the new S-SPADE algorithm which is in line to A-SPADE.
%The document also discusses the issue of complex variables.\\%[.5ex]

%\noindent
%\emph{Acknowledgment:}
%Thanks to Jakub Kudela for discussing the Augmented Lagrangian for complex variables.

\section{Alternating Direction Method of Multipliers (ADMM)}
\label{sec:ADMM}

%\noindent
We start with the brief review of ADMM \cite{Boyd2011ADMM}.
ADMM is able to solve problems of the form
\begin{equation}
\min_\x f(\x) + g(A\x),
\label{eq:problem_formulation1}
\end{equation}
where \(\x\in\CC^N\) and \(A:\CC^N\to \CC^P\) is a linear operator.
We assume $f,g$ real convex functions of (possibly complex) variables.
Problem \eqref{eq:problem_formulation1} can be reformulated as follows:
\begin{equation}
\min_{\x,\z} f(\x) + g(\z) \hspace{1em} \text{s.t.} \hspace{1em} A\x - \z = 0.
\label{eq:problem_formulation2}
\end{equation}
To solve optimization task \eqref{eq:problem_formulation2} we form the Augmented Lagrangian as
\begin{equation}
L_\rho(\x,\y,\z) = f(\x) + g(\z) + \y^{\top}(A\x-\z) + \frac{\rho}{2}\|A\x - \z\|^2_2,
\label{eq:augmented_lagrangian}
\end{equation}
where \(\rho > 0\) is called the \emph{penalty parameter}.
%\todo{sjednocovani indexu iteraci?}
%
ADMM consists of three steps:
\begin{subequations}\label{eq:ADMM_unscaled}
\begin{align}
\x^{(i+1)} &= \argmin_\x L_\rho \left(\x, \z^{(i)}, \y^{(i)}\right) \\
\z^{(i+1)} &= \argmin_\z L_\rho \left(\x^{(i+1)}, \z, \y^{(i)} \right) \\
\y^{(i+1)} &= \y^{(i)} + \rho \left( A\x^{(i+1)} - \z^{(i+1)}\right).
\end{align}
\end{subequations}

It is possible to convert ADMM to the so-called \emph{scaled form}, which is often more convenient;
this is done by defining the residual \(\rr = A\x - \z\).
In such a case, the last two terms of the Augmented Lagrangian in \eqref{eq:augmented_lagrangian} can be rewritten as
\begin{equation}
\y^\top\rr + \frac{\rho}{2}\|\rr\|^2_2 = \frac{\rho}{2}\|\rr + \frac{1}{\rho}\y\|^2_2 - \frac{1}{2\rho}\|\y\|^2_2 = 
\frac{\rho}{2}\|\rr+\u\|^2_2 - \frac{\rho}{2}\|\u\|^2_2,
\label{eq:scaled_form}
\end{equation}
where \(\u\) is a \emph{scaled dual variable} such that
%\(\u = \displaystyle\frac{1}{\rho}\y\)
\(\u = \y/\rho\)
(see Remark \ref{rem:scaled_form} for the proof of Eq.\,\eqref{eq:scaled_form}).

After the above manipulation, the Augmented Lagrangian in the scaled form is %can be written as
\begin{equation}
L_\rho(\x,\u,\z) = f(\x) + g(\z) + \frac{\rho}{2}\|\rr+\u\|^2_2 - \frac{\rho}{2}\|\u\|^2_2.
\label{eq:augmented_lagrangian_scaled}
\end{equation}
The scaled version of the ADMM is expressed as 
\begin{subequations}\label{eq:ADMM}
\begin{align}
\label{eq:ADMM_scaled_f}
\x^{(i+1)} & = \argmin_\x\left( f(\x) + \frac{\rho}{2} \|A\x - \z^{(i)} + \u^{(i)}\|^2_2\right) \\
\label{eq:ADMM_scaled_g}
\z^{(i+1)} & = \argmin_\z\left( g(\z) + \frac{\rho}{2} \|A\x^{(i+1)} - \z + \u^{(i)}\|^2_2\right) \\
\label{eq:ADMM_scaled_u}
\u^{(i+1)} & = \u^{(i)} + A\x^{(i+1)} - \z^{(i+1)}.
\end{align}
\end{subequations}

Note that in minimizing \(L_\rho(\x, \y, \z)\) over \(\x\) in the first step \eqref{eq:ADMM_scaled_f},
it was possible to omit \(g(\z)\). 
This term is not dependent on \(\x\), thus it does not play any role in finding the argument of the minima of \(L_\rho(\x, \y, \z)\).
The same argumentation applies to  \(\frac{\rho}{2}\|\u\|^2_2\).
%the other terms which were omitted.
Eq.\,\eqref{eq:ADMM_scaled_g} is obtained the same way.
%term \(\frac{\rho}{2}\|\u\|^2_2\) in \eqref{eq:augmented_lagrangian_scaled}.
%Similarly it is possible to omit terms \(f(\x)\) and \(\frac{\rho}{2}\|\u\|^2_2\) in \eqref{eq:ADMM_scaled_g}.

\begin{Rem}
	\label{rem:scaled_form}
	The second equality in Eq.\,\eqref{eq:scaled_form} is just a~substitution.
	The first one can be proven (right to left) using the following, assuming $\rr$ and $\y$ are real for the moment: 
	%\begin{subequations}
	\begin{equation}
		%\label{eq:scaled_form_proof}
		\begin{aligned}
		\label{eq:scaled_form_proof}
		\frac{\rho}{2}\|\rr+\frac{1}{\rho}\y\|^2_2-\frac{1}{2\rho}\|\y\|^2_2 &
		= \frac{\rho}{2}\langle\rr+\frac{1}{\rho}\y,\rr+\frac{1}{\rho}\y\rangle-\frac{1}{2\rho}\langle\y,\y\rangle \\&
		= \frac{\rho}{2}\left( \langle\rr,\rr\rangle+\langle\rr,\frac{1}{\rho}\y\rangle+\langle\frac{1}{\rho}\y,\rr\rangle+\langle\frac{1}{\rho}\y,\frac{1}{\rho}\y\rangle \right)-\frac{1}{2\rho}\langle\y,\y\rangle\\ &
		= \frac{\rho}{2}\|\rr\|^2_2+\frac{1}{2}\rr^\top\y+\frac{1}{2}\y^\top\rr\\ &
		= \frac{\rho}{2}\|\rr\|^2_2 + \y^\top\rr.
		\end{aligned}
		\end{equation}
	%\end{subequations}
	In the case of complex variables, %that the variables are assumed complex, 
	we need to guarantee the Augmented Lagrangian to be real (and therefore the minimization is legitimate). 
	For this reason, the Augmented Lagrangian in Eq.\,\eqref{eq:augmented_lagrangian}
	would need to be alternatively defined as
	\begin{equation}
	L_\rho(\x,\y,\z) = f(\x) + g(\z) + \y^{\top}
	\begin{bmatrix}
	\Re(A\x-\z)\\ \Im(A\x-\z)
	\end{bmatrix}
	+ \frac{\rho}{2}
	\norm{
	\begin{bmatrix}
	\Re(A\x-\z)\\ \Im(A\x-\z)
	\end{bmatrix}}_2^2,
	\label{eq:augmented_lagrangian_comlplex}
	\end{equation}
	now with $\y\in\RR^{2P}$.
	Then, \eqref{eq:scaled_form_proof} can be used again to prove Eq.\,\eqref{eq:scaled_form}, leading to Augmented Lagrangian in the form \eqref{eq:augmented_lagrangian_scaled}---the only difference is that $\rr,\u\in\RR^{2P}$ instead of $\CC^P$.
	In the implementation, this difference can be ignored and complex variables $\rr$ and $\u$ can be directly used since it holds
	\begin{equation}
	\norm{\c}_2^2 = \norm{\begin{bmatrix}
		\Re(\c) \\ \Im(\c)
		\end{bmatrix}}_2^2
	\end{equation}
	for any complex vector $\c$;
	therefore Equations \eqref{eq:ADMM_scaled_f} and \eqref{eq:ADMM_scaled_g}
	can be used as they are, and Eq.\,\eqref{eq:ADMM_scaled_u} corresponds only to addition of $\u$ and $\rr$,
	which can be equivalently performed either in $\CC^P$	or separately with the real and imaginary parts.
	%computing in $\RR^{2P}$.
\end{Rem}

\section{Sparse Audio Declipper (SPADE) -- Derivation using ADMM}
\label{sec:ASPADE}
SPADE algorithm(s) approximate the solution of the following non-convex, NP-hard synthesis- or analysis-regularized inverse problems
\cite{Kitic2015:Sparsity.cosparsity.declipping}:
\begin{equation}
	\label{eq:problem_ana}
	\min_{\x,\z} \norm{\z}_0  \quad \text{s.\,t.}
	\quad \x\in\Gamma(\y) \ \,\text{and}\ \,\norm{A\x-\z}_2\leq\epsilon \hspace{0.5em} \dots \hspace{0.5em} \text{A-SPADE}, 
\end{equation}
\begin{equation}
	\label{eq:problem_syn}
	\min_{\x,\z} \norm{\z}_0 \quad \text{s.\,t.}
	\quad \x\in\Gamma(\y) \ \,\text{and}\ \,\norm{\x-D\z}_2 \leq\epsilon \hspace{0.5em} \dots \hspace{0.5em} \text{S-SPADE},
\end{equation}
where \(\Gamma=\Gamma(\y)\) is the set of feasible solutions, \(\x\in\RR^N\) is signal in the time domain, \(\z\in\CC^P\) are signal coefficients.
Linear operator \(D: \CC^P \to \RR^N\) is the synthesis operator and \(A: \RR^N \to \CC^P\) is the analysis operator, and it holds \(D = A^*\).

Because of the computational reasons, we will restrict ourself exclusively to Parseval tight frames,
% as $D$ in the following,
i.e.\ it holds
\begin{equation}
	D\adjoint{D} = \adjoint{A}A = \Id.
	\label{eq:frame_operator}
\end{equation}
For such frames, it holds
\begin{equation}
	\label{eq:frame_inequality.synthesis}
	\norm{\adjoint{A}\c}_2 \leq \norm{\c}_2
	\quad \text{while}\quad 
	\norm{\adjoint{A}\c}_2=\norm{\c}_2 \ \text{for}\ \c\in\RA,
\end{equation}
where $\RA$ denotes the range space of $A$.
%and
%
%\begin{equation}
%	\label{eq:frame_inequality.analysis}
%	\norm{A\c}_2=\norm{\c}_2 \quad \text{for all}\ \c .
%\end{equation}
%\todo{this second relation is not utilized anywhere!?}

Problems \eqref{eq:problem_syn} and \eqref{eq:problem_ana} can also be written as a sum of two indicator functions, such that
\begin{equation}
\min_{\x, \z, k} \ig + \il \hspace{1em}\text{s.t.}\hspace{1em}
	\begin{array}{lll}
		 A\x - \z \hspace{-0.5em}&= 0 \hspace{0.5em} \dots \hspace{-0.3em} &\text{A-SPADE} \\
		 \x - D\z \hspace{-0.5em}&= 0 \hspace{0.5em} \dots \hspace{-0.3em} &\text{S-SPADE},
	\end{array}
\label{eq:problem_SPADE_ADMM}
\end{equation}
where \(\ig\) is indicator function forcing the result to lie in the set of feasible solutions~\(\Gamma\) 
and \(\il\) is a symbolic notation for indicator function \(\iota_{\{\z \hspace{0.3em}|\hspace{0.3em} \|\z\|_0 \leq k\}} (\z)\), that enforces the sparsity of the signal.

In the particular case of A-SPADE, our goal is to minimize the following %optimization problem
\begin{equation}
\argmin_{\x, \z, k} \underbrace{\ig}_{f(\x)} + \underbrace{\mystrut{4pt}\il}_{g(\z)} \hspace{1em}\text{s.t.}\hspace{1em}  A\x - \z = 0,
\label{eq:problem_ASPADE}
\end{equation}
which for fixed \(k\) correspond to problems of the form \eqref{eq:problem_formulation2}.
However, the difference is the non-convexity of $g$, therefore we aim only at an \emph{approximation} of \eqref{eq:problem_ASPADE} using ADMM \cite{Boyd2011ADMM}.
We form the Augmented Lagrangian %as
\begin{equation}
L_\rho(\x,\y,\z) = \ig + \il + \y^{\top}(A\x-\z) + \frac{\rho}{2}\|A\x - \z\|^2_2.
\label{eq:augmented_lagrangian2}
\end{equation}
Using the scaled form, the Augmented Lagrangian can be also written as
\begin{equation}
L_\rho(\x,\z, \u) = \ig + \il + \frac{\rho}{2}\|A\x - \z + \u \|^2_2 - \frac{\rho}{2}\|\u\|^2_2.
\label{eq:augmented_lagrangian3}
\end{equation}
According to definition of ADMM in \eqref{eq:ADMM}, we can form the ADMM algorithm for problem \eqref{eq:problem_ASPADE} as follows:
\begin{subequations}\label{eq:ADMM_ASPADE_unconst}
\begin{align}
\label{eq:ADMM_ASPADE_unconst_f}
\x^{(i+1)} &= \argmin_\x\left( \ig + \frac{\rho}{2} \|A\x - \z^{(i)} + \u^{(i)}\|^2_2\right) \\
\label{eq:ADMM_ASPADE_unconst_g}
\z^{(i+1)} &= \argmin_\z\left( \il + \frac{\rho}{2} \|A\x^{(i+1)} - \z + \u^{(i)}\|^2_2\right) \\
\label{eq:ADMM_ASPADE_unconts_u}
\u^{(i+1)} &= \u^{(i)} + A\x^{(i+1)} - \z^{(i+1)}.
\end{align}
\end{subequations}

Note that thanks to the indicator functions, the penalty parameter \(\rho\) does not play any role in finding the argument of the minima, therefore can be omitted. 
ADMM steps \eqref{eq:ADMM_ASPADE_unconst} can be also written in the constrained form, such that we get
\begin{subequations}\label{eq:ADMM_ASPADE_const}
\begin{align}
\label{eq:ADMM_ASPADE_const_f}
\x^{(i+1)} & = \argmin_\x\|A\x - \z^{(i)} + \u^{(i)}\|^2_2 \hspace{1em} \text{s.t.} \hspace{1em} \x \in \Gamma \\
\label{eq:ADMM_ASPADE_const_g}
\z^{(i+1)} & = \argmin_\z\|A\x^{(i+1)} - \z + \u^{(i)}\|^2_2 \hspace{1em} \text{s.t.} \hspace{1em} \|\z\|_0 \leq k\\
\label{eq:ADMM_ASPADE_conts_u}
\u^{(i+1)} & = \u^{(i)} + A\x^{(i+1)} - \z^{(i+1)}.
\end{align}
\end{subequations}
The \(\x\)-update, i.e.\ the steps \eqref{eq:ADMM_ASPADE_unconst_f} and \eqref{eq:ADMM_ASPADE_const_f},
is an orthogonal projection onto the set of feasible solutions \(\Gamma\).
The \(\z\)-update, steps \eqref{eq:ADMM_ASPADE_unconst_g} and \eqref{eq:ADMM_ASPADE_const_g},
is solved (not approximated!) by hard-thresholding operator \(\mathcal{H}_k\), which sets all but \(k\) largest components in the magnitude of the input vector to zero. This is implied by the property
\begin{equation}
	\argmin_{\z\in\CC^P,\ \norm{\z}_0\leq k}\norm{\z-\s}_2^2=\mathcal{H}_k(\s)
	%\quad\forall\s\in\CC^P,
	\label{eq:Hkisoptimal}
\end{equation}
for any fixed $\s\in\CC^P$.
\begin{Rem}
	It is straightforward to show the property \eqref{eq:Hkisoptimal}. First, it is clear that the optimal solution $\z_\mathrm{opt}$ has sparsity $k$, unless $\norm{\s}_0<k$ and then $\z_\text{opt}=\s$.
	Supposing that it is not the trivial case, we will divide $\norm{\z-\s}_2^2$ into two parts---the sum of squares of the coefficients of $\z$ that are forced to be zero by the condition $\norm{\z}_0\leq k$, and the rest.
	Apparently, the first part is minimized by setting the $P-k$ smallest coefficients of $\s$ to zero.
	Furthermore, the second part is minimized by setting the remaining $k$ coefficients
	equal to corresponding values of $\s$, leading us to the solution $\z_\mathrm{opt} = \mathcal{H}_k(\s)$.
\end{Rem}
\begin{Rem}
	Step \eqref{eq:ADMM_ASPADE_const_f} is a projection that seeks for a signal-domain solution.
	This is not trivial, however, along with \cite{Kitic2015:Sparsity.cosparsity.declipping},
 it can be transposed to a~more convenient form, which is
%The projection \eqref{eq:ADMM_ASPADE_const_f} is 
	\begin{equation}
		\x^{(i+1)} = \argmin_\x\|\x - \adjoint{A}(\z^{(i)} - \u^{(i)})\|^2_2 \hspace{1em} \text{s.t.} \hspace{1em} \x \in \Gamma,
		\tag{\ref{eq:ADMM_ASPADE_const_f}'}
		\label{eq:ADMM_ASPADE_const_f_2}
	\end{equation}	
	i.e.\ a~projection onto $\Gamma$ that is easy to implement.
	We will now prove the equivalence of \eqref{eq:ADMM_ASPADE_const_f} and \eqref{eq:ADMM_ASPADE_const_f_2}. % Tag 20a' is devided to 20a and '. Is it possible to have only 20a' refering to the equation 20a' ?
	
	Denote $\s = \z^{(i)} - \u^{(i)}$ and define its unique factorization
	\begin{equation}
	\s = \underbrace{A\xi}_{\in\RA} + \underbrace{\mystrut{2.1pt}\varepsilon.}_{\perp\RA}
	\label{eq:s_factor}
	\end{equation}
	Using this, Eq.\,\eqref{eq:ADMM_ASPADE_const_f} can be rewritten as
	\begin{equation}
		\x^{(i+1)} = \argmin_\x\|A(\x-\xi) - \varepsilon\|^2_2 \hspace{1em} \text{s.t.} \hspace{1em} \x \in \Gamma.
		\label{eq:ADMM_ASPADE_const_f_3}
	\end{equation}
	Because clearly $A(\x-\xi)\in\RA$ and $\varepsilon\perp\RA$, we can use the Pythagorean theorem and write
	\begin{equation}
	\x^{(i+1)} = \argmin_\x\|A(\x-\xi)\|^2_2 + \|\varepsilon\|^2_2\hspace{1em} \text{s.t.} \hspace{1em} \x \in \Gamma.
	\label{eq:ADMM_ASPADE_const_f_4}
	\end{equation}
	Since we are searching for the argument of the minima, we can omit the term $\|\varepsilon\|^2_2$, which is independent of $\x$.
	Next, we will make use of the property of Parseval tight frames \eqref{eq:frame_inequality.synthesis},
	%\begin{equation}
		%\norm{\c}=\norm{\adjoint{A}\c} \hspace{1em}\text{for}\hspace{1em}\c\in\RA.
	%\end{equation}
	leading to the equivalence of Eq.\,\eqref{eq:ADMM_ASPADE_const_f_4} and
	\begin{equation}
	\x^{(i+1)} = \argmin_\x\|\underbrace{\adjoint{A}A}_{\Id}(\x-\xi)\|^2_2\hspace{1em}\text{s.t.} \hspace{1em} \x \in \Gamma.
	\label{eq:ADMM_ASPADE_const_f_5}
	\end{equation}
	Finally, we will show that $\xi = \adjoint{A}\s$. Using $A\xi=\s-\varepsilon$
	%(this is clear from \eqref{eq:s_factor})
	and the property \eqref{eq:frame_operator}, we can write
	\begin{equation}
		\xi=\adjoint{A}(\s-\varepsilon)=\adjoint{A}\s-\adjoint{A}\varepsilon.
	\end{equation}
	We know that $\varepsilon\perp\RA\Leftrightarrow\langle\varepsilon,A\omega\rangle=0\;\forall\omega$. Using the definition of adjoint operator, $\langle\varepsilon,A\omega\rangle=\langle\adjoint{A}\varepsilon,\omega\rangle\;\forall\varepsilon\;\forall\omega$, therefore $\langle\adjoint{A}\varepsilon,\omega\rangle=0\;\forall\omega \Leftrightarrow\adjoint{A}\varepsilon=0$ and $\xi=\adjoint{A}\s$.
	
	This completes the proof, because then $\|{\x-\xi}\|^2_2 = \|{\x-\adjoint{A}\s}\|^2_2 = \|{\x-\adjoint{A}(\z^{(i)} - \u^{(i)})}\|^2_2$ and that leads to Eq.\,\eqref{eq:ADMM_ASPADE_const_f_2}.
\end{Rem}

The SPADE algorithms, as introduced in the original paper
\cite{Kitic2015:Sparsity.cosparsity.declipping},
are shown in Alg.\,\ref{alg:aspade} (A-SPADE) and Alg.\,\ref{alg:sspade} (S-SPADE), respectively.

\vspace{1em}
\noindent
\begin{minipage}{\linewidth}{
\small
\begin{minipage}[t]{0.49\linewidth}
\flushleft
%\scalebox{.99}{
\vspace{0pt}
\begin{algorithm}[H]
\DontPrintSemicolon
	\SetKwInput{KwRequire}{Require}
	\SetKw{KwReturn}{return}
	
	\KwRequire{\(A, \y, \Mr, \Mh, \Ml, s, r, \epsilon\)} \vspace{0.3em}
	\({\hat{\x}^{(0)} = \y, \u^{(0)}=\mathbf{0}, i=1, k=s}\)\;
	\(\bar{\z}^{(i+1)} = \mathcal{H}_k\left(A\hat{\x}^{(i)}+\u^{(i)}\right)\)\;
	\({\hat{\x}^{(i+1)} = \argmin_\x{\|A\x-\bar{\z}^{(i+1)}+\u^{(i)}\|_2^2}} \newline\text{s.t.\hspace{0.5em}}\x \in \Gamma\)\;
	\eIf{\(\|A\hat{\x}^{(i+1)}-\bar{\z}^{(i+1)}\|_2 \leq \epsilon\)}{\textup{terminate}\;}
	{\(\u^{(i+1)}=\u^{(i)}+A\hat{\x}^{(i+1)}-\bar{\z}^{(i+1)}\)\;
	\(i \leftarrow i+1\)\;
	\If{\(i\modulo r = 0\)}{\(k \leftarrow k+s\)\;}
	go to 2\;
	}
	\KwReturn{\(\hat{\x} = \hat{\x}^{(i+1)}\)}
	\caption{A-SPADE}
	\label{alg:aspade}
\end{algorithm}
%}
\end{minipage}
\hfill
\begin{minipage}[t]{0.49\textwidth}
\flushright
%\scalebox{.98}{
\vspace{0pt}
\begin{algorithm}[H]
\DontPrintSemicolon
	\SetKwInput{KwRequire}{Require}
	\SetKw{KwReturn}{return}
	
	\KwRequire{\(D, \y, \Mr, \Mh, \Ml, s, r, \epsilon\)} \vspace{0.3em}
	%\({\hat{\z}^{(0)} = D^{\textsf{H}}\y, \u^{(0)}=\mathbf{0}, i=1, k=s}\)\;
	\({\hat{\z}^{(0)} = D^*\y, \u^{(0)}=\mathbf{0}, i=1, k=s}\)\;
	\(\bar{\z}^{(i+1)} = \mathcal{H}_k\left(\hat{\z}^{(i)}+\u^{(i)}\right)\)\;
	\({\hat{\z}^{(i+1)} = \argmin_\z{\|\z-\bar{\z}^{(i+1)}+\u^{(i)}\|_2^2}}\newline\text{s.t.\hspace{0.5em}}	D\z \in \Gamma\)\;
	\eIf{\(\|\hat{\z}^{(i+1)}-\bar{\z}^{(i+1)}\|_2 \leq \epsilon\)}{\textup{terminate}\;}
	{\(\u^{(i+1)}=\u^{(i)}+\hat{\z}^{(i+1)}-\bar{\z}^{(i+1)}\)\;
	\(i \leftarrow i+1\)\;
	\If{\(i\modulo r = 0\)}{\(k \leftarrow k+s\)\;}
	go to 2\;
	}
	\KwReturn{\(\hat{\x} = D\hat{\z}^{(i+1)}\)}
	\caption{S-SPADE}
	\label{alg:sspade}
\end{algorithm}
%}
\end{minipage}
}
\end{minipage}

\vspace{1em}
\noindent
While the paper \cite{Kitic2015:Sparsity.cosparsity.declipping} does not present the derivation of either A-SPADE or S-SPADE, we have shown in Sections \ref{sec:ADMM} and \ref{sec:ASPADE} that 
A-SPADE perfectly fits the ADMM paradigm.
On the other hand, the S-SPADE version in Alg.\,\ref{alg:sspade} follows a~problem formulation which is different from the one the paper coped with.
However, it is still derived from ADMM;
it is easy to show that the problem formulation corresponding to the S-SPADE algorithm from \cite{Kitic2015:Sparsity.cosparsity.declipping} is
\begin{equation}
	\label{eq:problem_syn_orig}
	\min_{\w,\z} \norm{\z}_0 \quad \text{s.\,t.}
	\quad D\w\in\Gamma(\y) \ \,\text{and}\ \,\norm{\w-\z}_2 \leq\epsilon.
\end{equation}
This original S-SPADE can be derived from \eqref{eq:problem_syn_orig} in the same way we described the \mbox{A-SPADE} derivation, starting from  \eqref{eq:problem_ana}.
% For example, on could notice that the definition of ADMM implies that the two principal steps should minimize the Augmented Lagrangian over different variables. 
% The original S-SPADE keeps the Hard Thresholding step and only substitutes \(A\x\) with \(\z\) and therefore \(\x\) with \(D\z\).

Section \ref{sec:SSPADE-DR} introduces a new synthesis version of the SPADE algorithm %, which is correctly derived from the definition of ADMM (therefore the ``done right'' in the title).
derived from formulation \eqref{eq:problem_syn}, which is in some sense much more in line with the analysis approach than formulation \eqref{eq:problem_syn_orig}
(therefore the ``done right'' in the title).

\begin{Rem}
	Note that \cite{Kitic2015:Sparsity.cosparsity.declipping} reported that the projection step in S-SPADE is computationally expensive.
	Our recent paper \cite{ZaviskaRajmicPrusaVesely2018:RevisitingSSPADE} showed that there exists a one-step projection, making \mbox{A-SPADE} and S-SPADE identical from this point of view.
	It is a~kind of paradox that we were able to speed up an algorithm that is %not correctly derived,
	derived in an inconsistent way (with the initial formulation),
	although producing reasonable results.
	In the new S-SPADE presented below, our fast projection step is not helpful.
\end{Rem}

\section{Synthesis Version of SPADE, Done Right}
\label{sec:SSPADE-DR}

The problem formulation of this version of SPADE algorithm is very similar to \eqref{eq:problem_ASPADE}, only with the difference that \(\il\) will be identified with \(f(\z)\) and likewise \(\ig\) will be identified with \(g(\x)\). 
The primal variable, in this case, will be \(\z\) representing signal coefficients, opposite to A-SPADE, where the primal variable was signal waveform, \(\x\).
We formulate the problem as
\begin{equation}
\argmin_{\x, \z, k} \underbrace{\mystrut{4.2pt}\il}_{f(\z)} + \underbrace{\ig}_{g(\x)} \hspace{1em}\text{s.t.}\hspace{1em}  D\z - \x = 0,
\label{eq:problem_SSPADE_DR}
\end{equation}
We again form the Augmented Lagrangian for this problem as:
\begin{equation}
L_\rho(\x,\y,\z) = \il + \ig + \y^{\top}(D\z-\x) + \frac{\rho}{2}\|D\z - \x\|^2_2.
\label{eq:augmented_lagrangian_sspade1}
\end{equation}
As in \eqref{eq:augmented_lagrangian3}, we use the scaled form
\begin{equation}
L_\rho(\x,\z, \u) = \il + \ig + \frac{\rho}{2}\|D\z - \x + \u \|^2_2 - \frac{\rho}{2}\|\u\|^2_2,
\label{eq:augmented_lagrangian_sspade2}
\end{equation}
leading to the ADMM steps as follows:
\begin{subequations}\label{eq:ADMM_SSPADE_const}
\begin{align}
	\label{eq:ADMM_SSPADE_const_f}
	\z^{(i+1)} & = \argmin_\z\|D\z - \x^{(i)} + \u^{(i)}\|^2_2 \hspace{1em} \text{s.t.} \hspace{1em} \|\z\|_0 \leq k \\
	\label{eq:ADMM_SSPADE_const_g}
	\x^{(i+1)} & = \argmin_\x\|D\z^{(i+1)} - \x + \u^{(i)}\|^2_2 \hspace{1em} \text{s.t.} \hspace{1em} \x \in \Gamma\\
	\label{eq:ADMM_SSPADE_conts_u}
	\u^{(i+1)} & = \u^{(i)} + D\z^{(i+1)} - \x^{(i+1)}.
\end{align}
\end{subequations}
When comparing ADMM steps of A-SPADE and S-SPADE, one can notice a different order of the minimization steps.
This is caused by setting the variable \(\z\) as the primal variable in S-SPADE. 
According to \cite{Boyd2011ADMM}, it should be possible to freely choose the order of the \(f\)- and \(g\)-update steps.
As a consequence, it is possible to first apply the hard thresholding and then the projection onto \(\Gamma\) in both A-SPADE and S-SPADE,
being aware that the statement from \cite{Boyd2011ADMM} generally holds true for convex problems only.%

%Boyd also mentions several variations of ADMM.
%It is possible to perform several \(f\)-updates and \(g\)-updates before one update of the dual variable \(\u\). This makes the ADMM to be closer to the standard method of multipliers. 
%It is also possible to perform the dual update after both \(f\)- and \(g\)-updates.

The \(f\)-update \eqref{eq:ADMM_SSPADE_const_f} is a~challenging task to solve,
but we rely on the frequent ADMM behaviour that it still converges even when the individual steps
are computed (in contrast to A-SPADE) only approximately. 
%Therefore, we can rewrite the \eqref{eq:ADMM_SSPADE_const_f} as
%%
%\begin{equation}
%\z^{i+1} = \argmin_\z\|D^*D\z - D^*(\x^i - \u^i)\|^2_2 \hspace{1em} \text{s.t.} \hspace{1em} \|\z\|_0 \leq k, \\
%\label{eq:hard_thresholding_sspade}
%\end{equation}
%%
%since \(D^* = A\) does not change the norm.
%Here, \(D^*(\x^i - \u^i)\) is a fixed vector, which leads to sparse recovery (approximation) task.
%%
%And even knowing that the Gramm operator \(D^*D\) is not identity,
We interpret Eq.\,\eqref{eq:ADMM_SSPADE_const_f} as follows: we search for a $k$-sparse vector of coefficients $\z^{(i+1)}$, such that the synthesis operator $D$ applied to these coefficients gives signal closest to $(\x^{(i)}-\u^{(i)})$.
We suggest approximating the solution of this optimization task with the hard thresholding operator
\(\mathcal{H}_k\) applied to the analysis coefficients of $(\x^{(i)}-\u^{(i)})$, i.e.
\begin{equation}
%\z^{i+1} = \mathcal{H}_k \left(D^*(\x^i - \u^i)\right),
\z^{(i+1)} \approx \z^{(i+1)}_\mathrm{appr} = \mathcal{H}_k \left(D^*(\x^{(i)} - \u^{(i)})\right),
\label{eq:hard_thresholding_sspade2}
\end{equation}
and we will show that this approximation is sufficiently close to the proper solution of \eqref{eq:ADMM_SSPADE_const_f}.
Using the properties of Parseval tight frame 
\eqref{eq:frame_operator}
%\begin{equation}
	%\norm{D\x}_2\leq\norm{\x}_2
%\end{equation}
and \eqref{eq:frame_inequality.synthesis}, we can write
\begin{subequations}\label{eq:Hkisgood}
	\begin{align}
	\label{eq:Hkisgood_1}
	\norm{D\z^{(i+1)} - \x^{(i)} + \u^{(i)}}_2^2 &= \norm{D\z^{(i+1)} - D\adjoint{D}(\x^{(i)} - \u^{(i)})}_2^2 \\
	\label{eq:Hkisgood_2}
	&= \norm{D(\z^{(i+1)} - \adjoint{D}(\x^{(i)} - \u^{(i)}))}_2^2 \\
	\label{eq:Hkisgood_3}
	&\leq \norm{\z^{(i+1)} - \adjoint{D}(\x^{(i)} - \u^{(i)})}_2^2.
	\end{align}	
\end{subequations}
Using the property \eqref{eq:Hkisoptimal}, the norm \eqref{eq:Hkisgood_3} is minimized by $\z^{(i+1)}_\mathrm{appr}$ defined by \eqref{eq:hard_thresholding_sspade2}.
Therefore the error of the approximation in the time domain (i.e.\ the original norm in \eqref{eq:Hkisgood_1}) is bounded by the minimal value of \eqref{eq:Hkisgood_3} and thus we can expect it to be sufficiently small.

The solution to \eqref{eq:ADMM_SSPADE_const_g} is obtained by a simple projection in the time domain.

The final version of the S-SPADE done right is shown in Algorithm~\ref{alg:sspade_dr}.

Note that the computational complexity is identical in all three SPADEs
and it is dominated by the cost of the transforms.

\begin{Rem}
	Note that the ADMM step \eqref{eq:ADMM_SSPADE_const_f} can be also approximated using
	other means, for example by a greedy algorithm such as the OMP, where the signal to be sparsely approximated is $(\x^{(i)}-\u^{(i)})$, the dictionary is $D$ and the algorithm is forced to terminate after $k$~iterations.
\end{Rem}

\vspace{1em}
\noindent
\begin{minipage}{\linewidth}
	\small
	\centering
	\begin{minipage}[t]{0.57\linewidth}
		\vspace{0pt}
		\begin{algorithm}[H]
			\DontPrintSemicolon
			\SetKwInput{KwRequire}{Require}
			\SetKw{KwReturn}{return}
			
			\KwRequire{\(D, \y, \Mr, \Mh, \Ml, s, r, \epsilon\)} \vspace{0.3em}
			\({\hat{\x}^{(0)} = \y, \u^{(0)}=\mathbf{0}, i=0, k=s}\)\;
			\(\bar{\z}^{(i+1)} = \mathcal{H}_k\left(D^*(\hat{\x}^{(i)} - \u^{(i)})\right)\)\;
			\({\hat{\x}^{(i+1)} = \argmin_\x{\|D\bar{\z}^{(i+1)} -\x +\u^{(i)}\|_2^2}} \hspace{0.54em}\text{s.t.\hspace{0.5em}}\x \in \Gamma\hspace{-1em}\)\;
			\eIf{\(\|D\bar{\z}^{(i+1)}-\hat{\x}^{(i+1)}\|_2 \leq \epsilon\)}{\textup{terminate}\;}
			{\(\u^{(i+1)}=\u^{(i)}+D\bar{\z}^{(i+1)}-\hat{\x}^{(i+1)}\)\;
				\(i \leftarrow i+1\)\;
				\If{\(i\modulo r = 0\)}{\(k \leftarrow k+s\)\;}
				go to 2\;
			}
			\KwReturn{\(\hat{\x} = \hat{\x}^{(i+1)}\)}
			\caption{S-SPADE done right}
			\label{alg:sspade_dr}
		\end{algorithm}
		%}
	\end{minipage}
\end{minipage}
\vspace{1em}

\section{Equivalence of SPADEs in case of unitary operators}
Note also that the three presented SPADEs are equivalent for unitary operators $D = A^* = A^{-1}$.
This fact can be shown as an equivalence of formulations of the three problems---the resulting algorithms are equivalent as a~consequence,
 since the derivation process is the same for all three.
First, let us repeat the formulation \eqref{eq:problem_SPADE_ADMM} for A-SPADE and S-SPADE done right
\begin{equation}
\min_{\x, \z, k} \ig + \il \hspace{1em}\text{s.t.}\hspace{1em}
\begin{array}{lll}
A\x - \z \hspace{-0.5em}&= 0 \hspace{0.5em} \dots \hspace{-0.3em} &\text{A-SPADE} \\
\x - D\z \hspace{-0.5em}&= 0 \hspace{0.5em} \dots \hspace{-0.3em} &\text{S-SPADE},
\end{array}
\tag{\ref{eq:problem_SPADE_ADMM}}
\end{equation}
and let \eqref{eq:problem_syn_orig} be reformulated in similar form as
\begin{equation}
\min_{\w, \z, k} \iota_{\Gamma(\y)}(D\w) + \il \hspace{1em}\text{s.t.}\hspace{1em}
\w - \z = 0.
\label{eq:problem_SPADE_orig_ADMM}
\end{equation}
The two formulations in \eqref{eq:problem_SPADE_ADMM} differ only by the constraint binding $\x$ with $\z$
and in the unitary case, one can be easily translated into the other by applying $D = A^{-1}$ (A-SPADE to S-SPADE) or $A = D^{-1}$ (S-SPADE to A-SPADE) on both sides of the constraint.
To show the equivalence of \eqref{eq:problem_SPADE_orig_ADMM} and \eqref{eq:problem_SPADE_ADMM},
we substitute $\x = D\w$, leading to $\w = D^{-1}\x = A\x$.
Formulation \eqref{eq:problem_SPADE_orig_ADMM} then attains the form
\begin{equation}
\min_{\x, \z, k} \iota_{\Gamma(\y)}(\underbrace{DA}_{\Id}\x) + \il \hspace{1em}\text{s.t.}\hspace{1em}
A\x - \z = 0,
\end{equation}
which is the A-SPADE formulation from \eqref{eq:problem_SPADE_ADMM}.

%This document is intended as  reference
\pagebreak
{
\inputencoding{cp1250}
\bibliographystyle{IEEEtran}
\bibliography{literatura}      % Bibliography file (usually '*.bib' )
}

\end{document}